\documentclass{birkjour}
 \newtheorem{thm}{Theorem}[section]
 \newtheorem{cor}[thm]{Corollary}

 \theoremstyle{definition}
 
 \theoremstyle{remark}
 \newtheorem{rem}[thm]{Remark}
 
 \numberwithin{equation}{section}

\begin{document}

%
%
%
%
%
%
%
%
%

\title[Generalized Fibonacci Numbers]
 {Special Matrices Associated with Generalized Fibonacci Numbers}

\author[G. Cerda-Morales]{Gamaliel Cerda-Morales}
\address{Instituto de Matem\'aticas, Pontificia Universidad Cat\'olica de Valpara\'iso, Blanco Viel 596, Valpara\'iso, Chile.}
\email{gamaliel.cerda.m@mail.pucv.cl}



\date{January 1, 2004}


\begin{abstract}
In \cite{Ka}, the authors obtained a method for deriving special matrices, whose powers are related to Fibonacci and Lucas numbers. In the study, it has been developed a method for deriving special matrices of $3\times 3$ dimensions, whose powers are related to Horadam and generalized Fibonacci numbers, and some special matrices have been found via the method developed. 
\end{abstract}

\maketitle
\section{Introduction}
The Horadam numbers have many interesting properties and applications in many fields of science (see, e.g., \cite{Lar1}). The Horadam numbers $H_{n}(a,b;r,s)$ are defined by the recurrence relation
\begin{equation}\label{e1}
H_{0}=a,\ H_{1}=b,\ H_{n+2}=rH_{n+1}+sH_{n},\ n\geq0.
\end{equation}
Another important sequence is the generalized Fibonacci sequence $\{h_{n}\}_{n\in \mathbb{N}}$ defined by $h_{n}=H_{n}(0,1;r,s)$.  

Let $\mathbf{Q}_{h}=\left[
\begin{array}{cc}
r& s \\ 
1& 0  
\end{array}
\right]$ be a companion matrix of the generalized Fibonacci sequence $\{h_{n}\}_{n\geq0}$ defined by the two-order linear recurrence relation
\begin{equation}\label{eq:1}
h_{0}=0,\ h_{1}=1,\ \ \ h_{n}=r h_{n-1}+s h_{n-2}\ (n\geq2).
\end{equation}
Then, by an inductive argument, the $n$-th power $\mathbf{Q}_{h}^{n}$ has the matrix form
\begin{equation}\label{eq:2}
\mathbf{Q}_{h}^{n}=\left[
\begin{array}{ccc}
h_{n+1}& s h_{n} \\ 
h_{n}&s h_{n-1}
\end{array}
\right]\ (n\geq1).
\end{equation}
This property provides an alternate proof of the Cassini-type formula for $\{h_{n}\}_{n\geq0}$, using $\det(\mathbf{Q}_{h}^{n})=\det(\mathbf{Q}_{h})^{n}$:
\begin{equation}\label{eq:3}
h_{n}^{2}-h_{n-1}h_{n+1}=(-s)^{n-1}.
\end{equation}

It is well-known that the usual generalized Fibonacci numbers can be expressed using Binet's formula
\begin{equation}\label{eq:4}
h_{n}=\frac{\alpha_{r,s}^{n}-\beta_{r,s}^{n}}{\alpha_{r,s}-\beta_{r,s}}
\end{equation}
where $\alpha_{r,s}$ and $\beta_{r,s}$ are the roots of the equation $x^{2}-rx-s=0$. Furthermore, $2\alpha_{r,s}=r+A_{h}$ and $2\beta_{r,s}=r-A_{h}$, where $$A_{h}=\sqrt{r^{2}+4s} \in \mathbb{R}.$$ 

From the Binet formula in Eq. (\ref{eq:4}) and using the classic identities
\begin{equation}\label{e0}
\alpha_{r,s}+\beta_{r,s}=r,\ \ \alpha_{r,s}\beta_{r,s}=-s,
\end{equation}
we have for any integer $n\geq1$:
\begin{align*}
h_{n}-\beta_{r,s} h_{n-1}&=\frac{\alpha_{r,s}^{n}-\beta_{r,s}^{n}}{\alpha_{r,s}-\beta_{r,s}}-\frac{\beta_{r,s} \left(\alpha_{r,s}^{n-1}-\beta_{r,s}^{n-1}\right)}{\alpha_{r,s}-\beta_{r,s}}\\
&=\frac{\alpha_{r,s}^{n-1}(\alpha_{r,s}-\beta_{r,s})}{\alpha_{r,s}-\beta_{r,s}}\\
&=\alpha_{r,s}^{n-1}.
\end{align*}
Then, we obtain 
\begin{equation}\label{eq:5}
h_{n}-\beta_{r,s} h_{n-1}=\alpha_{r,s}^{n-1}\ (n\geq1).
\end{equation}
Multipling Eq. (\ref{eq:5}) by $\alpha_{r,s}$, using $\alpha_{r,s}\beta_{r,s}=-s$, and if we change $\alpha_{r,s}$ and $\beta_{r,s}$ role above process, we obtain the lineal approximation of $\{h_{n}\}_{n\geq0}$
\begin{equation}\label{eq:6}\textrm{Lineal app. of $\{h_{n}\}$}: \left\{
\begin{array}{c }
\alpha_{r,s}^{n}=\alpha_{r,s} h_{n}+s h_{n-1}\\
\beta_{r,s}^{n}=\beta_{r,s} h_{n}+s h_{n-1}.
\end{array}
\right.
\end{equation}

In Eq. (\ref{eq:6}), if we change $\alpha_{r,s}$ and $\beta_{r,s}$ into the companion matrix $\mathbf{Q}_{h}$ and change $h_{n-1}$ into the matrix $h_{n-1}\mathbf{I}_{2}$, where $\mathbf{I}_{2}$ is the $2\times 2$ identity matrix, then we obtain the matrix form Eq. (\ref{eq:2}) of $\mathbf{Q}_{h}^{n}$
\begin{equation}
\mathbf{Q}_{h}^{n}=h_{n}\mathbf{Q}_{h}+s h_{n-1}\mathbf{I}_{2}\left(=\left[
\begin{array}{ccc}
h_{n+1}& s h_{n} \\ 
h_{n}&s h_{n-1}
\end{array}
\right]\right).
\end{equation}

It is known that some properties of Horadam and generalized Fibonacci numbers can be proved by using matrices (see, e.g., \cite{Cer}).  This study has two-stages: First, an approach to the derivation of $3\times 3$ dimensional matrix whose powers are related to generalized Fibonacci numbers is presented. Then, based on this approach, some special $3\times 3$ dimensional matrices are obtained and some related identities are given.

\section{Some Special Matrices}

In this section, by giving an approach for derivating $3\times 3$ dimensional matrices associated with generalized Fibonacci numbers, based on this approach, special matrices are obtained.

\subsection{Methods and Results}
Suppose that $\mathbf{A}_{h}=\left[
\begin{array}{ccc}
a_{11}& a_{12}& a_{13} \\ 
a_{21}&a_{22}& a_{23} \\ 
a_{31}&a_{32}& a_{33} 
\end{array}
\right]$ is any $3\times 3$ matrix such that the eigenvalues of it are $\alpha_{r,s}$, $\beta_{r,s}$ and $0$, and all the entries of it are integers. Represent the eigenvectors corresponding the eigenvalues  $\alpha_{r,s}$, $\beta_{r,s}$ and $0$, respectively as $$\mathbf{x}=\left[
\begin{array}{c}
x_{1} \\ 
x_{2}\\ 
x_{3}\end{array}
\right],\ \mathbf{y}=\left[
\begin{array}{c}
y_{1} \\ 
y_{2}\\ 
y_{3}\end{array}
\right],\ \mathbf{z}=\left[
\begin{array}{c}
z_{1} \\ 
z_{2}\\ 
z_{3}\end{array}
\right].$$
For the eigenvalue-eigenvector pairs $(\alpha_{r,s},\mathbf{x})$, $(\beta_{r,s},\mathbf{y})$ and $(0,\mathbf{z})$, respectively, it is obtained the linear equations systems $\mathbf{A}_{h}\mathbf{x}=\alpha_{r,s}\mathbf{x}$,  $\mathbf{A}_{h}\mathbf{y}=\beta_{r,s}\mathbf{y}$ and  $\mathbf{A}_{h}\mathbf{z}=0\mathbf{z}$.

The matrix $\mathbf{A}_{h}$ is diagonalizable since it has different eigenvalues. In this case, without loss of the generality, it can be written $\mathbf{A}_{h}=\mathbf{P}\mathbf{D}_{h}\mathbf{P}^{-1}$, where $$\mathbf{D}_{h}=\left[
\begin{array}{ccc}
\alpha_{r,s}& 0& 0 \\ 
0&\beta_{r,s}& 0 \\ 
0&0&0\\
\end{array}
\right]$$ and and $\mathbf{P}$ is a nonsingular matrix. From this, we get $\mathbf{A}_{h}^{n}=\mathbf{P}\mathbf{D}_{h}^{n}\mathbf{P}^{-1}$ for all integers $n\geq 1$. Then, $$\mathbf{A}_{h}^{n}=\mathbf{P}\left[
\begin{array}{ccc}
\alpha_{r,s}^{n}& 0& 0 \\ 
0&\beta_{r,s}^{n}& 0 \\ 
0&0&0\\
\end{array}
\right]\mathbf{P}^{-1}.$$
Considering Eq. (\ref{eq:6}), it is obtained
\begin{align*}
\mathbf{P}^{-1}\mathbf{A}_{h}^{n}\mathbf{P}
&=h_{n}\left[
\begin{array}{ccc}
\alpha_{r,s}& 0& 0 \\ 
0&\beta_{r,s}& 0 \\ 
0&0&0\\
\end{array}
\right]+sh_{n-1}\left[
\begin{array}{ccc}
1& 0& 0 \\ 
0&1& 0 \\ 
0&0&1
\end{array}
\right]-sh_{n-1}\left[
\begin{array}{ccc}
0& 0& 0 \\ 
0&0& 0 \\ 
0& 0 &1\end{array}
\right]\\
&=h_{n}\mathbf{D}_{h}+sh_{n-1}\mathbf{I}_{3} -sh_{n-1}\left[
\begin{array}{cccc}
0& 0& 0 \\ 
0&0& 0\\ 
0& 0 &1\end{array}
\right],
\end{align*}
that is,
\begin{align*}
\mathbf{A}_{h}^{n}&=h_{n}\mathbf{P}\mathbf{D}_{h}\mathbf{P}^{-1}+sh_{n-1}\mathbf{P}\mathbf{I}_{3}\mathbf{P}^{-1}-sh_{n-1}\mathbf{P}\left[
\begin{array}{ccc}
0& 0& 0\\ 
0&0& 0 \\ 
0& 0 &1\end{array}
\right]\mathbf{P}^{-1}\\
&=h_{n}\mathbf{A}_{h}+s h_{n-1}\mathbf{I}_{3}-sh_{n-1}\mathbf{P}\left[
\begin{array}{ccc}
0& 0& 0 \\ 
0&0& 0 \\ 
0&0&1\end{array}
\right]\mathbf{P}^{-1},
\end{align*}
where $\mathbf{I}_{3}$ is the $3\times 3$ identity matrix. 

Thus, it is seen that the power of the matrix $\mathbf{A}_{h}$ is associated with generalized Fibonacci numbers. Now, if we write $$\mathbf{E}=\mathbf{P}\left[
\begin{array}{ccc}
0& 0& 0 \\ 
0&0& 0\\ 
0& 0 &1\end{array}
\right]\mathbf{P}^{-1}=\left[
\begin{array}{ccc}
0& 0&z_{1} \\ 
0&0&z_{2}\\ 
0&0&z_{3}\\
\end{array}
\right]\mathbf{P}^{-1},$$ then we get
\begin{equation}\label{ee1}
\mathbf{A}_{h}^{n}=h_{n}\mathbf{A}_{h}+sh_{n-1}\left(\mathbf{I}_{3}-\mathbf{E}\right).
\end{equation}

\subsection{Particular Cases for generalized Fibonacci Numbers}

Now, it is given some special matrices $\mathbf{A}_{h}$ which occur for the special cases of the eigenvectors $\mathbf{x}$, $\mathbf{y}$ and $\mathbf{z}$ and some identities associated with these. If it is chosen $$\mathbf{x}_{1}=\left[
\begin{array}{c}
x_{1} \\ 
x_{2}\\ 
x_{3}\end{array}
\right]=\left[
\begin{array}{c}
\alpha_{r,s} \\ 
\beta_{r,s}\\ 
-1\end{array}
\right],\ \mathbf{y}_{1}=\left[
\begin{array}{c}
y_{1} \\ 
y_{2}\\ 
y_{3}
\end{array}
\right]=\left[
\begin{array}{c}
\beta_{r,s}\\
\alpha_{r,s} \\  
-1\end{array}
\right],$$
then, it is necessary to hold the systems $\mathbf{A}_{h}\mathbf{x}_{1}=\alpha_{r,s}\mathbf{x}_{1}$ and $\mathbf{A}_{h}\mathbf{y}_{1}=\beta_{r,s}\mathbf{y}_{1}$ for being these vectors are eigenvectors associated with the eigenvalues $\alpha_{r,s}$ and $\beta_{r,s}$, respectively, of the matrix $\mathbf{A}_{h}$. From these, the following equations must be hold:
\begin{equation}\label{z1}
\left\lbrace \begin{aligned} \alpha_{r,s} a_{11}+\beta_{r,s} a_{12}-a_{13}&=\alpha_{r,s}^{2}\\
\beta_{r,s} a_{11}+\alpha_{r,s} a_{12}-a_{13}&=\beta_{r,s}^{2}\end{aligned}\right . ,
\end{equation}
\begin{equation}\label{z2}
\left\lbrace \begin{aligned} \alpha_{r,s} a_{21}+\beta_{r,s} a_{22}-a_{23}&=-s\\
\beta_{r,s} a_{21}+\alpha_{r,s} a_{22}-a_{23}&=-s\end{aligned}\right . 
\end{equation}
and
\begin{equation}\label{z3}
\left\lbrace \begin{aligned} \alpha_{r,s} a_{31}+\beta_{r,s} a_{32}-a_{33}&=-\alpha_{r,s}\\
\beta_{r,s} a_{31}+\alpha_{r,s} a_{32}-a_{33}&=-\beta_{r,s}\end{aligned}\right . .
\end{equation}
Or equivalent to
\begin{equation}\label{total}
\left\lbrace \begin{aligned} ra_{11}+ra_{12}-2a_{13}=r^{2}+2s&,\ a_{11}-a_{12}=r\\
ra_{21}+ra_{22}-2a_{23}=-2s&,\ a_{21}-a_{22}=0\\
ra_{31}+ra_{32}-2a_{33}=-r&,\ a_{31}-a_{32}=-1\end{aligned}\right . .
\end{equation}

Now, let us consider the different choices of the vector $\mathbf{z}$. For example, if it is chosen as $z_{1}=z_{3}=t$ and $z_{2}=-t$, where $t\in \mathbb{Z}^{\times}$ is arbitrary, matrix $\mathbf{E}$ in Eq. (\ref{ee1}) is obtained as
\begin{align*}
\mathbf{E}&=\left[
\begin{array}{ccc}
0& 0&t \\ 
0&0&-t\\ 
0&0&t\\
\end{array}
\right]\left[
\begin{array}{ccc}
 \alpha_{r,s} &  \beta_{r,s} &t \\ 
 \beta_{r,s} & \alpha_{r,s} &-t\\ 
-1&-1&t\\
\end{array}
\right]^{-1}\\
&=\frac{1}{r} \left[
\begin{array}{ccc}
1& 1&r \\ 
-1&-1&-r\\ 
1&1&r\\
\end{array}
\right],
\end{align*}
with $r\in \mathbb{R}^{\times}$, and also, the equation $\mathbf{A}_{h}\mathbf{z}=0\mathbf{z}$ turns to the system 
\begin{equation}\label{z5}
\left\lbrace \begin{aligned} a_{11}-a_{12}+a_{13}&=0\\
a_{21}-a_{22}+a_{23}&=0\\
a_{31}-a_{32}+a_{33}&=0\\\end{aligned}\right . .
\end{equation}

From the solutions of the equations systems (\ref{total}) and (\ref{z5}), we get 
\begin{equation}\label{z10}
\mathbf{A}_{h}=\frac{1}{r} \left[
\begin{array}{ccc}
r(r-1)+s& s-r &-r^{2} \\ 
-s&-s&0\\ 
1-r&1&r\\
\end{array}
\right]
\end{equation}

If we use the equality (\ref{ee1}) for the matrix $\mathbf{A}_{h}$ in (\ref{z10}), then we obtain
\begin{align*}
\mathbf{A}_{h}^{n}&=h_{n}\mathbf{A}_{h}+sh_{n-1}\left(\mathbf{I}_{3}-\mathbf{E}\right)\\
&=\frac{h_{n}}{r} \left[
\begin{array}{ccc}
r(r-1)+s& s-r &-r^{2} \\ 
-s&-s&0\\ 
1-r&1&r\\
\end{array}
\right]\\
&\ \ +\frac{sh_{n-1}}{r} \left[
\begin{array}{ccc}
r-1& -1&-r \\ 
1&r+1&r\\ 
-1&-1&0\\
\end{array}
\right]\\
&=\frac{1}{r}\left[
\begin{array}{ccc}
h_{n+2}-h_{n+1}&-(h_{n+1}-sh_{n})&-rh_{n+1} \\ 
-s(h_{n}-h_{n-1})&s(h_{n-1}-sh_{n-2})&rsh_{n-1}\\ 
-(h_{n+1}-h_{n})&h_{n}-sh_{n-1}&rh_{n}\\
\end{array}
\right],
\end{align*}
or all integers $n\geq 1$. Thus, we have been proved the following theorem:

\begin{thm}
If $\mathbf{A}_{h}=\frac{1}{r} \left[
\begin{array}{ccc}
r(r-1)+s& s-r &-r^{2} \\ 
-s&-s&0\\ 
1-r&1&r\\
\end{array}
\right]$. Then,
\begin{equation}\label{res1}
\mathbf{A}_{h}^{n}=\frac{1}{r}\left[
\begin{array}{ccc}
h_{n+2}-h_{n+1}&-(h_{n+1}-sh_{n})&-rh_{n+1} \\ 
-s(h_{n}-h_{n-1})&s(h_{n-1}-sh_{n-2})&rsh_{n-1}\\ 
-(h_{n+1}-h_{n})&h_{n}-sh_{n-1}&rh_{n}\\
\end{array}
\right],\ n\geq 1.
\end{equation}
\end{thm}

\begin{rem}
Since one of the eigenvalue of the matrix $\mathbf{A}_{h}$ is zero, the matrix $\mathbf{A}_{h}$ is singular. So, the result obtained is valid for only all integers $n\geq 1$ and $h_{-1}=1/s$ for convenience.
\end{rem}

\begin{cor}
For all $n\geq 1$. The following equalities are satisfied
\begin{equation}\label{g1}
\left[
\begin{array}{ccc}
1& 0&-1 \\ 
-1&-1&0\\ 
0&1&1\\
\end{array}
\right]^{n}=\left[
\begin{array}{ccc}
F_{n}&-F_{n-1}&-F_{n+1} \\ 
-F_{n-2}&F_{n-3}&F_{n-1}\\ 
-F_{n-1}&F_{n-2}&F_{n}\\
\end{array}
\right],
\end{equation}
\begin{equation}\label{g2}
\left(\frac{1}{2} \left[
\begin{array}{ccc}
3& -1 &-4 \\ 
-1&-1&0\\ 
0&1&2\\
\end{array}
\right]\right)^{n}=\frac{1}{2}\left[
\begin{array}{ccc}
P_{n+2}-P_{n+1}&-(P_{n+1}-P_{n})&-2P_{n+1} \\ 
-(P_{n}-P_{n-1})&P_{n-1}-P_{n-2}&2P_{n-1}\\ 
-(P_{n+1}-P_{n})&P_{n}-P_{n-1}&P_{n}\\
\end{array}
\right],
\end{equation}
\begin{equation}\label{g3}
\left[
\begin{array}{ccc}
2& 1&-1 \\ 
-2&-2&0\\ 
0&1&1\\
\end{array}
\right]^{n}=\left[
\begin{array}{ccc}
2J_{n}&-(J_{n+1}-2J_{n})&-J_{n+1} \\ 
-4J_{n-2}&2(J_{n-1}-2J_{n-2})&2J_{n-1}\\ 
-2J_{n-1}&J_{n}-2J_{n-1}&J_{n}\\
\end{array}
\right],
\end{equation}
where $F_{n}=h_{n}(0,1;1,1)$, $P_{n}=h_{n}(0,1;2,1)$ and $J_{n}=h_{n}(0,1;1,2)$ is the $n$-th Fibonacci, Pell and Jacobsthal number, respectively.
\end{cor}

Taking different choices of the eigenvectors $\mathbf{x}$, $\mathbf{y}$ and $\mathbf{z}$, and progressing similarly to the above, it can be given different matrices $\mathbf{A}_{h}$ and related results. Now, it is presented, without proof, some of these kind of results, for clarification.

For the choice of $(\alpha_{r,s},\left[\begin{array}{c}\alpha_{r,s}\\ \beta_{r,s}\\ -1\end{array}\right])$, $(\beta_{r,s},\left[\begin{array}{c}\beta_{r,s}\\ \alpha_{r,s}\\ -1\end{array}\right])$ and $(0,\left[\begin{array}{c}t \\t \\ -t\end{array}\right])$:

\begin{thm}
If $\mathbf{A}_{h}=\frac{1}{r-2} \left[
\begin{array}{ccc}
r(r-1)+s& r+s &r^{2}+2s \\ 
-s&-s&-2s\\ 
1-r&-1&-r\\
\end{array}
\right]$. Then,
\begin{equation}\label{res2}
\mathbf{A}_{h}^{n}=\frac{1}{r-2}\left[
\begin{array}{ccc}
h_{n+2}-h_{n+1}&h_{n+1}+sh_{n}&h_{n+2}+sh_{n} \\ 
-s(h_{n}-h_{n-1})&-s(h_{n-1}+sh_{n-2})&-s(h_{n}+sh_{n-2})\\ 
-(h_{n+1}-h_{n})&-(h_{n}+sh_{n-1})&-(h_{n+1}+sh_{n-1})\\
\end{array}
\right],
\end{equation}
where $r\in \mathbb{R}$ and $r\notin \{0,2\}$.
\end{thm}

For the choice of $(\alpha_{r,s},\left[\begin{array}{c}\alpha_{r,s}\\ \beta_{r,s}\\ -1\end{array}\right])$, $(\beta_{r,s},\left[\begin{array}{c}\beta_{r,s}\\ \alpha_{r,s}\\ -1\end{array}\right])$ and $(0,\left[\begin{array}{c}-t \\t \\ t\end{array}\right])$:

\begin{thm}
If $\mathbf{A}_{h}=\frac{1}{r} \left[
\begin{array}{ccc}
r(r+1)+s& r+s &r^{2} \\ 
-s&-s&0\\ 
-(r+1)&-1&-r\\
\end{array}
\right]$. Then,
\begin{equation}\label{res3}
\mathbf{A}_{h}^{n}=\frac{1}{r}\left[
\begin{array}{ccc}
h_{n+2}+h_{n+1}&h_{n+1}+sh_{n}&rh_{n+1} \\ 
-s(h_{n}+h_{n-1})&-s(h_{n-1}+sh_{n-2})&-rsh_{n-1}\\ 
-(h_{n+1}+h_{n})&-(h_{n}+sh_{n-1})&-rh_{n}\\
\end{array}
\right].
\end{equation}
\end{thm}

Taking determinant in Eq. (\ref{res3}) and using that $\det(\mathbf{A}_{h})=0$, we obtain
\begin{cor}
For $n\geq 2$, we have
\begin{equation}
h_{n}^{3}+h_{n-1}^{2}h_{n+2}+h_{n+1}^{2}h_{n-2}=h_{n}(h_{n-2}h_{n+2}+2h_{n-1}h_{n+1}).
\end{equation}
\end{cor}

\section{Summary}
Based on the approach given in \cite{Ka}, it is seen that it can be written any finite dimensional matrix and related results. Although the approach is simple, we believe that it is important in terms of its role and useful in the study related to these subjects. For example, using the generalized sequence $h_{n}(0,1;r,s)$ we can derive a matrix approximation of a large number of sequences, such as Fibonacci, Pell, Jacobsthal and Balancing. In a future study we intend to apply the method to recursive sequences of order 3; in particular, the Tribonacci and third-order Jacobsthal numbers.


\end{document}